\documentclass[12pt]{amsart}
\usepackage{amsmath,amsthm,amsfonts,amssymb,latexsym,mathabx,url, bm}

\usepackage{mathrsfs}
\usepackage{subdepth} 

\newtheorem{theorem}{Theorem}[section]
\newtheorem{prop}[theorem]{Proposition}
\newtheorem{lem}[theorem]{Lemma}

\newtheorem{cor}[theorem]{Corollary}

\numberwithin{equation}{section}

\makeatletter
\renewcommand{\pmod}[1]{\allowbreak\mkern7mu({\operator@font mod}\,\,#1)}
\makeatother

\newcommand{\cE}{\mathcal{E}}
\newcommand{\cL}{\mathcal{L}}
\newcommand{\N}{\mathbb{N}}
\newcommand{\sP}{\mathscr{P}}

\newcommand{\T}{{\mathscr{T}}}

\newcommand{\Vol}{\operatorname{Vol}}
\newcommand{\meas}{\operatorname{meas}}

\renewcommand{\ssum}[1]{\sum_{\substack{#1}}}  %%% stacked sum
\newcommand{\eps}{\ensuremath{\varepsilon}}

\renewcommand{\le}{\leqslant}
\renewcommand{\leq}{\leqslant}
\renewcommand{\ge}{\geqslant}
\renewcommand{\geq}{\geqslant}

\newcommand{\fl}[1]{{\ensuremath{\left\lfloor {#1} \right\rfloor}}}  % floor
\newcommand{\cl}[1]{{\ensuremath{\left\lceil #1 \right\rceil}}} % ceiling
\newcommand{\order}{\asymp}      %% order of magnitude

\newcommand{\pfrac}[2]{\left(\frac{#1}{#2}\right)}  %%% frac with auto-sized paren

\renewcommand{\bm}{\ensuremath{\mathbf{m}}}

\newcommand{\bxi}{\ensuremath{\boldsymbol{\xi}}}
\newcommand{\cM}{\mathcal{M}}
\newcommand{\tY}{\widetilde{Y}}
\newcommand{\ttheta}{\widetilde{\theta}}
\newcommand{\tr}{\widetilde{r}}
\newcommand{\tH}{\widetilde{H}}

\newcommand{\RR}{{\mathbb R}}

\begin{document}

\title{Integers with many divisors in a dyadic interval}
\author{Kevin Ford, Mikhail R. Gabdullin, Kyle Pratt}

\address{Department of Mathematics, University of Illinois at Urbana-Champaign} 
\email{ford126@illinois.edu}

\address{Department of Mathematics, University of Illinois at Urbana-Champaign} 
\email{mikhail.gabdullin@gmail.com}

\address{Department of Mathematics, Brigham Young University}
\email{kyle.pratt@math.byu.edu}

\begin{abstract}
We find the correct order of magnitude for the number of $n\leq x$ having at least $r$ divisors in the interval $(y,2y]$, uniformly in a wide range of the parameters $x,y,r$.
\end{abstract}		
 
\date{\today}	
\maketitle

\section{Introduction}
Let $\tau(n;y,z)$ be the number of divisors of $n$ in the interval $(y,z]$.
For $3\leq y\leq z\leq x$ and $r\in\N$, we set
\begin{align*}
H_r(x,y,z) &= \# \big\{ n\le x: \tau(n;y,z)=r \big\},   \\
H_{\ge r}(x,y,z) &= \#\big\{n\leq x: \tau(n;y,z)\ge r \big\}.
\end{align*}
Define the constant
\[
\cE  =1-\frac{1+\log\log 2}{\log 2}=0.08607...
\]
In this paper we'll concentrate on the central case $z=2y$.

In \cite{F1}, the first named author determined the order of $H(x,y,z)=H_{\ge 1}(x,y,z)$
for all $x,y,z$, in particular showing that
\[
H_{\ge 1}(x,y,2y) \order \frac{x}{(\log Y)^{\cE}(\log\log Y)^{3/2}}, \quad Y=\min(3y,3x/y).
\]
It is also proved in \cite{F1} that for each fixed $r$, $H_r(x,y,2y) \order_r H_{\ge 1}(x,y,2y)$ when $y\le x^{1/2-\eps}$ and $y>x^{1/2+\eps}$ ($\eps>0$ fixed),
there being subtle issues when $y$ is close to $\sqrt{x}$.  The bounds have an unspecified
dependence on $r$.

Our goal in this paper is to bound $H_{\ge r}(x,y,2y)$ uniformly in a large range of $r$.
Define
\begin{align*}
E(x) = \pfrac{1+x}{\log 2}\log\pfrac{e\log 2}{1+x}, \qquad
\eta & =E'(0) = \frac{\log\log 2}{\log 2}=-0.5287...	
\end{align*}

\begin{theorem}\label{th1.1}
Fix $\delta>0$.  Let $100\leq y\leq x/100$, and assume that $Y=\min(y,x/(2y))$ is sufficiently
large.  Let $1\leq r\le (\log Y)^{\log 4-1-\delta}$ and define $\theta$ by $r=(\log Y)^{\theta}$.  Then
\begin{equation*}
H_{\ge r}(x,y,2y) \order_\delta \frac{x \log (2r)}{(\log Y)^{1-E(\theta)}(\log\log Y)^{3/2}},	
\end{equation*}	
the implied constants in $\order$ depending only on $\delta$. 
\end{theorem}

In the special case of smaller $r$, namely $1\le r\le \exp \{ C(\log\log Y)^{1/2} \},$
we have
\[
E(\theta)=E(0)+\theta E'(0)+O_C(1/\log\log Y)=1-\cE+\eta \theta + O_C(1/\log\log Y),
\]
and therefore
\[
H_{\ge r}(x,y,2y) \order_{C} \frac{x r^{\eta} \log(2r)}{(\log Y)^{\cE}(\log\log Y)^{3/2}},
\]
the lower bound holding if $Y$ is sufficiently large in terms of $C$.

\subsection{Heuristic argument}
We adapt the heuristic argument from \cite{F1,F3}.  Assume $10\le y\le \sqrt{x}$.
%For brevity, let $\tau(n,y,z)$ be the number of divisors of $n$ in $(y,z]$.
Write $n=n'n''$, where $n'$ is composed only of primes $\le 2y$ and
$n''$ is composed only of primes $>2y$.
For simplicity, assume $n'$ is squarefree and $n'\le y^{100}$.
Assume for the moment that the set $D(n')=\{ \log d : d|n' \}$ is
uniformly distributed in $[0,\log n']$.  If $n'$ has $k$ prime factors, then 
the expected value of $\tau(n';y,2y)$ should be about 
$\frac{2^k\log 2}{\log n'} \asymp
\frac{2^k}{\log y}$.  This is $\ge r$ precisely when $k\ge k_0 + O(1)$,
where $k_0 := \fl{\frac{\log\log y+\log r}{\log 2}}$.
Using the fact (e.g. Theorem 08 of \cite{Divisors}) 
that the number of $n\le x$ with $n'$ having $k$ prime factors is
of order
$$
\frac{x}{\log y} \frac{(\log\log y)^k}{k!},
$$
we obtain a heuristic estimate for $H_{\ge r}(x,y,2y)$ of order
$$
\frac{x}{\log y} \sum_{k\ge k_0+O(1)} \frac{(\log\log y)^k}{k!} 
\asymp \frac{x (\log\log y)^{k_0}}{k_0! \log y} 
\asymp \frac{x (\log y)^{E(\theta)-1}} {(\log\log y)^{1/2}}, \quad \theta=\frac{\log r}{\log \log y}.
$$

This is slightly too big, and the reason stems from the uniformity
 assumption about $D(n')$.  Similar to the heuristics from \cite{F1,F3}, $D(n')$
 will be roughly uniformly distributed only for $n'$ with
 a special anatomy, namely those integers where $D(n')$ doesn't have
 any dense clusters of primes in intervals of length $\gg 1$.
 The probability that $n'$ has this property is about $\frac{\log r}{\log\log y}$. 
 
The above heuristic is fine for $r < (\log y)^{\log 4-1-\delta}$, but the
method breaks down when $r$ is near $(\log y)^{\log 4-1}$ due to the effect of integers
with $k<k_0$ prime factors.  Our upper bounds for the contribution of such
integers, Lemma \ref{lem:Lr}(b) below, are too weak to handle these integers.
It is possible that there is a phase transition of $H_{\ge r}(x,y,2y)$
in the vicinity of $(\log y)^{\log 4-1}$.  There are similarities between
the behavior of $H_{\ge r}(x,y,2y)$ and the behavior of $H(x,y,2^{1/r}y)$, which is also
dominated by integers with about $k_0$ prime factors below $y$.  The latter is
known to undergo a phase transition near $r=(\log y)^{\log 4-1}$
(see \cite[Chapter 2]{Divisors} for the basic theory, and \cite{FT} for a 
refined analysis for $H_{\ge 2}(x,y,2^{1/r}y)$).

\subsection{Notation.}
We adopt the following notational conventions: $\tau(n)$ is the number of positive divisors
of $n$, $\omega(n)$ is the number of distinct prime factors of $n$,
$P^-(n)$ is the smallest prime factor of $n$ with $P^-(1):=\infty$,
$P^+(n)$ is the largest prime factor of $n$ with $P^+(1):=0$.
The notations $f=O(g)$ and $f\ll g$ mean that $|f|/g$ is bounded.  When
$f,g$ are functions of more than one variable, dependence of the implied constants
on some variables is indicated by subscripts to $O,\ll$.
The symbol $p$, with or without subscripts, always denotes a prime.
For real $x$, $\fl{x}$ is the largest integer $\le x$ (the floor function) and $\cl{x}$ is the smallest
integer $\ge x$ (the ceiling function).

\subsection{The primary reduction}

\medskip
For real $r>0$ we define
\[
\cL_{\ge r}(a)=\{x\in\RR: \tau(a;e^x,2e^x)\geq r\}, \quad L_{\ge r}(a)=\meas\cL_{\ge r}(a).
\]
Note that 
\begin{equation}\label{1.1}
\cL_{\ge r}(a)=\cL_{\ge \cl{r}}(a)
\end{equation} 
for any $r>0$; in particular, $L_{\ge r}(a)=L_{\ge 1}(a)$ when $0<r\leq 1$. For $z\geq2$, define the set
\begin{equation}\label{1.2}
\sP(z)=\{a\in\N: P^+(a)\leq z, a\text{ squarefree}\}.
\end{equation}
Also let
\[
\vartheta(r,z) := \frac{\log r}{\log\log z}, \qquad \text{ so that } r=(\log z)^{\vartheta(r,z)}.
\]

Theorem \ref{th1.1} can be reduced to the following statement about sums of $L_{\ge t}(a)$.
It is a kind of global-to-local principle.

\begin{prop}\label{prop1.2}
\begin{enumerate}
\item[(a)] Fix $\delta'>0$.  For $z\ge 3$ and $1\leq t \leq (\log z)^{\log 4-1-\delta'}$, 
\[
\sum_{a\in \sP(z)}\frac{L_{\ge t}(a)}{a} \ll_{\delta'} \frac{(\log z)^{1+E(\vartheta(t,z))}\log (2t)}{(\log\log z)^{3/2}}.
\]
\item[(b)] There is a $z_0$ so that uniformly for $z\ge z_0$ and $1\leq t \leq (\log z)^{10}$, 
\[
\sum_{a\le z}\frac{L_{\ge t}(a)}{a} \gg \frac{(\log z)^{1+E(\vartheta(t,z))}\log (2t)}{(\log\log z)^{3/2}}.
\]
\end{enumerate}
\end{prop}

In Section 2 we deduce Theorem \ref{th1.1} from Proposition \ref{prop1.2}, and in Section 3 we prove Proposition \ref{prop1.2}. Throughout, all the implied constants in $O$- and $\ll-$symbols
 are absolute unless indicated by a subscript, e.g. $\ll_h$.

\section{Preliminaries} 

We start with several auxiliary statements. First, we will often use the standard estimates for the rough numbers (see, for instance, \cite[Theorem III.6.4]{ten}): one has
\begin{equation}\label{Phi-lower}
\#\{3x/4<n\leq x: P^-(n)>z\} \gg \frac{x}{\log z}	
\end{equation}
uniformly in $x\geq 2z$ and $z$ sufficiently large, and
\begin{equation}\label{Phi-upper}
\#\{n\leq x: P^-(n)>z\} \ll \frac{x}{\log z}
\end{equation}
for all $x\ge z\ge 2$.  Note that the condition $P^-(n)>z$ includes $n=1$.

The Prime Number Theorem with classical error term \cite[Chapter 18]{Da} implies that 
for all $A>1$ and $z\geq 3$, 
\begin{equation}\label{PNT}
\sum_{p\leq z}\frac1p=\log\log z+c_0+O_A\bigg(\frac{1}{(\log z)^{A}}\bigg),
\end{equation}
where $c_0$ is a constant.

Now we list some upper bounds for the $L_{\ge r}(a)$ function.

\newpage

\begin{lem}\label{lem:Lr}
Let $r\geq1$ be integer. We have
\begin{itemize}
\item[(a)] $L_{\ge r}(a)\leq \min\big\{\frac{3\tau(a)}{r}, \log (2a)\big\}$;	

\smallskip 

\item[(b)] $L_{\ge r}(ab) \leq \tau(b)L_{\ge r/\tau(b)}(a)$ whenever $\gcd(a,b)=1$; 	

\smallskip 

\item[(c)] For any primes $p_1<\cdots <p_k$,
\[
L_{\ge r}(p_1\cdots p_k) \leq  2^{k+2}\min\left(\frac1r,\min_{0\leq j\leq k} 2^{-j}\log(2p_1\cdots p_j) \right) .
\]
\end{itemize}	
\end{lem}

\begin{proof} a) First of all, $x\in \cL_{\ge r}(a)$ trivially implies  $-\log 2\le x\leq \log a$. Further, we have $x\in \cL_{\ge r}(a)$ if and only if there are divisors $d_1,\ldots,d_r$ of $a$ with $e^x<d_1<\cdots<d_r\leq 2e^x$. Hence
\[
\cL_{\ge r}(a)=\bigcup_{\substack{d_1<\cdots<d_r < 2d_1 \\ d_i|a\; (1\le i\le r)}}[\log d_r-\log 2, \log d_1),
\] 
which is a (potentially empty) union of half-open intervals $[x,x+\alpha_x)$ of length $\alpha_x\leq \log2$ (since $\log d_1\leq \log d_r$).  It follows that the open intervals
$I_x=(x-0.1,x+\log 2+0.1)$ form a covering of $\cL_{\ge r}(a)$.  The Vitali covering lemma
\cite[Lemma 3.15]{Folland} then guarantees the existence of a disjoint subcollection $I_{x_1},\ldots,I_{x_k}$ such that
$$
L_{\ge r}(a)\leq 3\sum_{j=1}^k|I_{x_j}|\leq 3(\log2+0.2)k \le 3k. 
$$
By construction, each of the $k$ intervals $(e^{x_j-0.1},e^{x_j+\log 2 +0.1})$ contains at least $r$ divisors of $a$. Since these intervals are mutually disjoint, we have $k\leq \tau(a)/r$. This concludes the proof of the part (a).

%Letting $\e$ go to zero, we conclude the proof of the part (a). 
	
\smallskip 
	
b) If $x\in \cL_{\ge r}(ab)$, then there are divisors $d_1,\ldots,d_r$ of $ab$ with 
\[
e^x<d_1<\cdots<d_r\leq 2e^x.
\] 	
For each $i$, we can write $d_i=d_i'd_i''$, where $d_i'|a, d_i''|b$ and $\gcd(d_i',d_i'')=1$. By the pigeonhole principle, there is a $d|b$ with $d=d_i''$ for at least $r/\tau(b)$ indices $i$. The corresponding (distinct) divisors $d_i'$ of $a$ then lie in the interval $(e^x/d,2e^x/d]$, and thus $x-\log d\in \cL_{\ge r/\tau(b)}(a)$. Therefore,
\[
\cL_{\ge r}(ab) \subseteq \bigcup_{d|b}\big\{y+\log d: y\in  \cL_{\ge r/\tau(b)}(a)\big\},
\]	
and the claim follows from \eqref{1.1}.
		
\smallskip 
	
(c) By part (a), $L_{\ge r}(p_1\cdots p_k) \le 2^{k+1}/r$.  
We may also invoke (b).  For each $j$, write $p_1\cdots p_k=ab$ with $a=p_1\cdots p_j$ 
and $b=p_{j+1}\cdots p_k$, and apply (b) plus
$L_{\ge r/\tau(b)}(a) \le \log(2a)$ which comes from (a).
\end{proof}	

%We note that, in Lemma 3(c), one can also get the more general estimate 
%$$
%L_{\ge r}(p_1\cdots p_k) \ll  \min_{0\leq j\leq k} 2^{k-j}\min\left(2^k/r,\log^+ %(2^{k-j}p_1\cdots p_j/r) \right),
%$$
%but we will not need it.  

%\smallskip 

For a nonnegative arithmetic function $f$ and $k\in\N$, define the function
\[
f_k(n)=\max\{f(sn): \gcd(s,n)=1, \, s\text{ squarefree}, \, \omega(s)\leq k\}.
\]
In particular, $f(n)\leq f_k(n)$ for any $k$. The next lemma is a variant of \cite[Lemma 2.2]{Kou17}.  

\begin{lem}\label{lem:suma}
For any fixed $h\in\N$, any nonnegative arithmetic function $f$ and $u,z\ge 2$,
\[
\sum_{a\in \sP(u)}\frac{f(a)}{a\log^h(P^+(a)+z/a)} \ll_h \frac{1}{\log^h z}\sum_{a\in \sP(u)} \frac{f_{h+2}(a)}{a}.
\]
\end{lem}

\begin{proof}
Let $w=\min(u,z^{1/4})$ and $\sP_1=\{a\in \sP(w): a>z^{1/2}\}$. The sum in the question is 	
\begin{equation}\label{sum-Pu}
\leq \sum_{a\in \sP_1}\frac{f(a)}{a\log^h P^+(a)}+\frac{4^h}{\log^h z}\sum_{a\in \sP(u)}\frac{f(a)}{a}.
\end{equation}
We write every $a\in\sP_1$ as $a=pb$, where $p=P^+(a)$; such $p$ exists because $a>z^{1/2}>1$. 
Since $a>z^{1/2}$ and $p\le z^{1/4}$, $b>z^{1/4}$ and hence
\begin{align}\label{sum-P1}
\sum_{a\in \sP_1}\frac{f(a)}{a\log^hP^+(a)} \leq &\sum_{p\leq w}\frac{1}{p\log^h p}
\ssum{b\in \sP(p)\\ b>z^{1/4}}\frac{f_1(b)}{b}\\
\leq & \frac{4^{h+1}}{(\log z)^{h+1}}\sum_{p\leq w}\frac{1}{p\log^h p}
\sum_{b\in \sP(p)}\frac{f_1(b)\log^{h+1}b}{b}. \nonumber 
\end{align}
Further,
\begin{align*}
\sum_{b\in \sP(p)}\frac{f_1(b)\log^{h+1}b}{b} &= \sum_{b\in \sP(p)}\frac{f_1(b)}{b}\sum_{p_1,\ldots,p_{h+1}|b}(\log p_1)\cdots(\log p_{h+1})\\ & \leq \sum_{p_1,\ldots,p_{h+1}\leq p}\frac{(\log p_1)\cdots(\log p_{h+1})}{\text{lcm}[p_1,\ldots,p_{h+1}]}\sum_{t\in \sP(p)}\frac{f_{h+2}(t)}{t},
\end{align*}
where we wrote $b=\text{lcm}[p_1,\ldots,p_{h+1}]t$. Since
\[
\sum_{p_1,\ldots,p_{h+1}\leq p}\frac{(\log p_1)\cdots(\log p_{h+1})}{\text{lcm}[p_1,\ldots,p_{h+1}]} 
\ll_h \log^{h+1}p, 
\]
from \eqref{sum-P1} we get
\begin{align*}
\sum_{a\in \sP_1}\frac{f(a)}{a\log^hP^+(a)} &\ll_h \frac{1}{(\log z)^{h+1}}\sum_{p\le w}
\frac{\log p}{p} \sum_{t\in \sP(p)} \frac{f_{h+2}(t)}{t}\\
&\ll \frac{\log w}{(\log z)^{h+1}}  \sum_{t\in \sP(u)} \frac{f_{h+2}(t)}{t}
\le \frac{1}{(\log z)^h}  \sum_{t\in \sP(u)} \frac{f_{h+2}(t)}{t}.
\end{align*}
Together with \eqref{sum-Pu}, this completes the proof.
\end{proof}

\begin{cor}\label{cor:sum-Lta-over-logs}
For $t\ge 1$, $u\ge 2$, $z\ge 2$ we have
\[
\sum_{a\in\sP(u)} \frac{L_{\ge t}(a)}{a\log^2 (P^+(a)+z/a)} \ll \frac{1}{\log^2 z} \sum_{a\in\sP(u)}
\frac{L_{\ge t/16}(a)}{a}.
\]
\end{cor}

\begin{proof}
Let $f(a)=L_{\ge t}(a)$ and $h=2$.  By Lemma \ref{lem:Lr} (b), $f_4(a)\le 16L_{\ge t/16}(a)$.
Now apply Lemma \ref{lem:suma}.
\end{proof}

Let $\omega(n,a,b)$ denote the number of prime divisors of $n$ in the interval $(a,b]$.

\begin{lem}\label{lem:l-prime-factors}
Let $l\in \N$ and set $\eta=\exp(-l/10)$. If $z\ge 2$ and $X>z^{\eta}$ then
\[
\#\{n\leq X: P^-(n)>z^\eta, \omega(n,z^\eta,z)\geq l-1\} 
\ll \frac{X}{\log z}\,e^{-l}.
\]  	
\end{lem} 

\begin{proof}
Let $l_0$ be a large constant.
When $l\le l_0$, \eqref{Phi-upper} implies that the left side is 
$\ll X/\log(z^{\eta})\ll e^{-l} X/\log z$.
When $l>l_0$, apply Theorem 3 of \cite{F4} with $T_1$ the (possibly empty) set of primes $\le z^\eta$
and $T_2$ the set of primes in $(z^\eta,\min(z,X)]$.  Let 
\[
H = \sum_{z^\eta<p\le \min(z,X)} \frac{1}p, \qquad H' = \sum_{p\le z^{\eta}} \frac{1}{p}.
\]
Then \eqref{PNT} implies $H\le l/10+O(1) \le \frac{l-1}{2}-2$ if $l_0$ is large enough. 
 Hence, the left side in the lemma is
\[
\ll e^{-H-H'} X \sum_{k\ge l-1} \frac{(H+2)^k}{k!}
\ll \frac{e^{l/10-H} X}{\log z} \cdot \frac{(H+2)^{l-1}}{(l-1)!} \ll \frac{e^{-l}X}{\log z}
\]
using Stirling's formula in the last step.
\end{proof}

\section{Deducing Theorem \ref{th1.1} from Proposition \ref{prop1.2}}

We begin with the lower bound in Theorem \ref{th1.1}. 
We assume that $Y$ is large enough.
Define
\[
\hat{\cL}(a)=\{x\in\RR: \tau(a;e^x,2^{1/2}e^x)\ge r\}, \quad \hat{L}(a)=\meas \, \hat{\cL}(a).
\]

We consider the numbers $n=apqb \in (3x/4,x]$ with $p,q$ prime and
\[
a\leq Y^{1/4}<p<q\leq \frac{Y^{3/4}}{2} < P^-(b), \qquad \log\pfrac{Y}{pq}\in \hat{\cL}(a).
\]
The last condition implies that $\tau(n;Y,2^{1/2} Y) \ge \tau(pqa;Y,2^{1/2}Y)\geq r$.
If $Y=\frac{x}{2y}$ and $d|n$ with $Y<d\le 2^{1/2} Y$, we have
$\frac{3x}{4d} < \frac{n}{d} \le \frac{x}{d}$ and hence $(3/\sqrt{8})y<\frac{n}{d} < 2y$.
Thus, in either case ($Y=y$ or $Y=\frac{x}{2y}$) we have $\tau(n;y,2y)\ge r$, as desired.

Trivially, $Y/(pq)\geq \frac12$. This gives
$$
\frac{x}{apq}\geq \frac{x}{2Y^{5/4}} \geq Y^{3/4},
$$
where we used that $Y^2\leq x/2$. Thus, for fixed $a,p,q$ as above, the number of possible $b$ is $\gg \frac{x}{apq\log Y}$ by \eqref{Phi-lower} and therefore
\begin{equation}\label{Hr-lower-1}
H_{\ge r}(x,y,2y) \gg \frac{x}{\log Y}\sum_{a\leq Y^{1/4}}\frac1a \;
\ssum{Y^{1/4}<p<q\leq (1/2)Y^{3/4}\\ \log(Y/(pq))\in \hat{\cL}(a)} \frac{1}{pq}.
\end{equation}
The set $\hat{\cL}(a)$ is a union of at most $\tau(a)$ disjoint intervals
$(v,v+w]$ with $w\le \frac{\log 2}{2}$. For each such interval 
$[v,v+w)$, we have $v\geq-\frac{\log2}{2}$ and $v+w\leq \log a\leq \frac14 \log Y$ and thus
\[
\ssum{Y^{1/4}<p<q\leq \frac12 Y^{3/4}\\ e^v \leq Y/(pq)<e^{v+w}}\; \frac{1}{pq}\geq 
\sum_{4Y^{1/4}\leq p\leq Y^{3/8}}\;\frac1p\;\sum_{Ye^{-v-w}/p<q\leq Ye^{-v}/p}\;\frac1q.
\]
By \eqref{PNT} and that $Ye^{-v-w}/p\le Y^{3/4}$, the inner sum on $q$ is
\[
\log\bigg(1+\frac{w}{\log(Ye^{-v-w}/p)}\bigg) + O\pfrac{1}{\log^{100}Y} \ge 
\frac{w}{\log Y} - O\pfrac{1}{\log^{100}Y}.
\]
Applying \eqref{PNT} once more, the sum on $p$ is $\ge \log \frac32-O(1/\log Y)\ge 0.2$
for large $Y$.  Summing over all intervals $[v,v+w)$ we conclude that
\[ 
\ssum{Y^{1/4}<p<q\leq (1/2)Y^{3/4}\\ \log(Y/(pq))\in \hat{\cL}(a)} \frac{1}{pq}
\geq \frac{0.2 \hat{L}(a)}{\log Y}-O\bigg(\frac{1}{(\log Y)^{100}}\bigg).
\]
Since $x\in \cL_{\ge 2r}(a)$ implies $x\in\hat{\cL}(a)$ or $x+\frac{\log 2}{2}\in\hat{\cL}(a)$,
\[
\hat{L}(a) \ge \frac{L_{\ge 2r}(a)}{2}.
\]
Therefore,
\[
H_{\ge r}(x,y,2y) \geq \frac{0.1x}{\log^2 Y}\sum_{a\leq Y^{1/4}}\frac{L_{\ge 2r}(a)}{a}-O\left(\frac{x}{(\log Y)^{100}}\sum_{a\leq Y^{1/4}}\frac{\tau(a)}{a}\right).
\]
Since $\sum_{a\leq Y^{1/4}}\tau(a)/a \ll \log^2 Y$, the lower bound in Theorem \ref{th1.1} now follows from Proposition \ref{prop1.2}(b) with $t=2r$ and $z=Y^{1/4}$, 
since  $2r \le 2(\log Y)^{\log 4-1} \le (\frac14\log Y)^{2}$ for large $Y$
and the right side in Theorem \ref{th1.1} is $\gg x/(\log Y)^{2}$.

\bigskip 

To establish the upper bound in Theorem \ref{th1.1}, we may assume without loss of
generality that $Y$ is sufficiently large in terms of $\delta$.
We need to  count  squarefree numbers with many divisors in a dyadic interval,
and begin with two preliminary sums.

\begin{lem}\label{lem:suml1}
Fix $\delta''>0$.  For positive integers $l$ let $z_l=z^{\exp(-l/10)}$.
Uniformly for $z$ sufficiently large in terms of $\delta''$
 and $1\le r\le (\log z)^{\log4 - 1-\delta''}$, we have
\begin{align}
\sum_{l=1}^r e^{-l} \; \sum_{a\in\sP(z)} \frac{L_{\ge r/(160l^2)}(a)}{a} &\ll_{\delta''}
\frac{\log(2r)(\log z)^{1+E(\vartheta(r,z))}}{(\log\log z)^{3/2}}, \label{suml-1}\\
\ssum{1\le l\le r \\ z_l\ge 2}  \; \sum_{a\in\sP(z_l)} \frac{L_{\ge r/(160l^2)}(a)}{a} &\ll_{\delta''}
\frac{\log(2r)(\log z)^{1+E(\vartheta(r,z))}}{(\log\log z)^{3/2}}. \label{suml-2}
\end{align}
\end{lem}

\begin{proof}
We may assume that $\delta''<1/100$.
We begin with \eqref{suml-1}.  Let $\theta''=\vartheta(r,z)$.
Fix $1\le l\le r$, let $t=\cl{r/(160l^2)}$ and $\theta'=\vartheta(t,z)$.
Then $t\le r\le (\log z)^{\log 4-1-\delta''}$ and $\theta'\le \theta''\le \log4-1-\delta''$.
Thus the  hypotheses of Proposition \ref{prop1.2} (b) hold.
Since $E(x)$ is decreasing,
\[
E(\theta')-E(\theta'') \ll (\theta''-\theta') \ll\; \frac{1+\log l}{\log\log z}
\]
and hence $(\log z)^{E(\theta')}\le (\log z)^{E(\theta'')} (2l)^{O(1)}$.  
By Proposition \ref{prop1.2} (a) with $\delta'=\delta''$, the sum in \eqref{suml-1} is
\[
\ll_{\delta'} \frac{(\log(2r))(\log z)^{1+E(\theta'')}}{(\log\log z)^{3/2}} \sum_{l=1}^r
e^{-l} (2l)^{O(1)} \ll  \frac{(\log(2r))(\log z)^{1+E(\theta'')}}{(\log\log z)^{3/2}}.
\]

Next, we bound the sum in \eqref{suml-2}.  Let $\delta'=\frac12 \delta''$.
Note that $z_l\ge 2$ is equivalent  to 
\[
l \le l^* := 10(\log\log z-\log\log 2).
\]
Let $\eps>0$ be sufficiently small in terms of $\delta''$ and set
$u_l=z^{\exp(-\eps l)}$.
Replace $\sP(z_l)$ in the sum in \eqref{suml-2} with $\sP(u_l)$.
With $l\le l^*$ fixed, set
\[
t=\cl{r/(160l^2)}, \qquad \theta'=\vartheta(t,u_l).
\]
Then 
\[
t \le r \le (\log z)^{\log 4-1-\delta''} \le (\log u_l)^{\log 4-1-\delta'}
\]
if $\eps$ is small enough.  We apply Proposition \ref{prop1.2} (a), obtaining
\begin{equation}\label{zlul}
\sum_{a\in\sP(z_l)} \frac{L_{\ge t}(a)}{a} \le 
\sum_{a\in\sP(u_l)} \frac{L_{\ge t}(a)}{a} \ll_{\delta'} \frac{(\log u_l)^{1+E(\theta')}\log(2r)}{(\log\log z)^{3/2}}.
\end{equation}
Now
\[
\theta' = \frac{\log r + O(1+\log l)}{\log\log z-\eps l} \ge \theta'' + O\pfrac{1+\log l}{\log\log z},
\]
Since $E(x)$ is decreasing,
\[
E(\theta') \le E(\theta'') + O\pfrac{1+\log l}{\log\log z}.
\]
Since $\theta'\le \log 4-1$,  $E(\theta')+1\ge 3-\log 4\ge 1.6$.  It follows that
\begin{align*}
(\log u_l)^{1+E(\theta')} &\le e^{-1.6\eps l}(\log z)^{1+E(\theta')} \ll
(2l)^{O(1)} e^{-1.6\eps l} (\log z)^{1+E(\theta'')}.
\end{align*}
Inserting this into \eqref{zlul} and summing on $l$ yields \eqref{suml-2}.
\end{proof}

\begin{lem}\label{lem:H-tilde-short-int}
Define
\[
\tH_{\ge r}(x,y,z)=\#\big\{n\le x: \tau(n;y,2y)\ge r, n\text{ is squarefree}  \big\}.
\]
Fix $\delta''>0$.
For $Y$ large enough in terms of $\delta''$, with $Y=\min(y,x/(2y))$,
and for
\[
r \leq (\log Y)^{\log 4-1-\delta''},
\]
we have
\[
\tH_{\ge r}(x,y,2y)-\tH_{\ge r}(x/2,y,2y) \ll_{\delta''} \frac{x \log(2r)}
{(\log Y)^{1-E(\vartheta(r,Y))}(\log\log Y)^{3/2}}.
\]
\end{lem}

\begin{proof}
The number of integers $n\le x$ divisible by $2y$ (if $2y$ is an integer) is $\le x/(2y)\le x/Y$,
which is acceptable since $E(\vartheta(r,Y))\in(0,1)$.

For squarefree $n$ counted in the left side not divisible by $2y$, let $y<d_1<\cdots <d_r< 2y$ be the $r$ smallest divisors
of $n$ in $(y,2y)$, and let $e_i=n/d_i$, so that $x/(4y) < e_r<\cdots< e_1<x/y$. There are two cases: 
\begin{enumerate}
\item[(i)] $P^+(e_i)<P^+(n)$ for at least $r/2$ divisors $e_i$; 
\item[(ii)] $P^+(d_i)<P^+(n)$ for at least $r/2$ divisors $d_i$.
\end{enumerate} 

In case (i) let $z=\frac{x}{4y}$ and $\cM = \{e_i : P^+(e_i)<P^+(n) \}$, while in case (ii)
let $z=y$ and $\cM =  \{d_i : P^+(d_i)<P^+(n) \}$.  In each case,
\[
z \ge \frac{Y}{2}, \qquad |\cM| \ge r/2, \qquad \forall\, m\in \cM, z<m<4z.
\]
The remaining argument is identical in both cases, so we'll just fix $z$ and $\cM$ as above.
Let 
\[
\{p_1, p_2,\ldots, p_s\} = \{P^+(m) : m\in \cM \}, \quad p_1>p_2>\cdots > p_s.
\]
For $1\le l\le s$, let $S_l = \{ m\in \cM: P^+(m)=p_l \}.$  Since 
$\sum_l |S_l| = |\cM| \ge r/2$, there is some $l$ with 
\[
|S_l| \ge \frac{r}{5l^2}.
\]
Fix $l$, let $\beta=1/(160 l^2)$ and factor 
\[
n=ap_l b, \qquad P^+(a) < p_l < P^-(b).
\]
By assumption, $p_l<P^+(n)$, so $b>p_l$.  For each $m\in \cM$, $(m/p_l)|a$.
Thus, since $z<m<4z$, 
\begin{equation}\label{apz}
\text{either } \log(z/p_l) \text{ or } \log(2z/p_l) \text{ lies in } \cL_{\ge 16\beta r}(a).
\end{equation}
Also $p_l \ge m/a > z/a$ and hence
\begin{equation}\label{plQ}
p_l > Q(a) := \max \big( P^+(a), z/a \big).
\end{equation}
For fixed $a$ and $p_l$ satisfying \eqref{apz} and \eqref{plQ},
\begin{equation}
p_l < b \le \frac{x}{ap_l}, \quad b\text{ squarefree}, \quad \omega(b,p_l,z)\ge l-1, \quad P^-(b)>p_l.
\end{equation}
If we ignore the condition $\omega(b,p_l,z)\ge l-1$, \eqref{Phi-upper}
implies that the number of possibilities for $b$ is 
\begin{equation}\label{b-upper-1}
\ll \frac{x}{ap_l \log Q(a)}.
\end{equation}
On the other hand, if $p_l \ge z^{\exp\{-l/10\}}$, then Lemma \ref{lem:l-prime-factors} implies that the number of possible $b$ is
\[
\ll \frac{e^{-l}x}{ap_l \log z}.
\]
For a given $a$, $\cL_{\ge 16\beta r}(a)$ is the union of at most $\tau(a)$ disjoint 
intervals of length $\le \log 2$.  Applying \eqref{PNT} to each shows that
\[
\sum_{c=1}^2 \ssum{p>Q(a) \\ \log(cz/p)\in \cL_{\ge 16\beta r}(a)} \frac{1}{p} \ll
\frac{L_{\ge 16\beta r}(a)}{\log Q(a)} + \frac{\tau(a)}{(\log Q(a))^{100}}.
\]
If follows from Corollary \ref{cor:sum-Lta-over-logs} (with $u=z$ and with
$u=z^{\exp(-l/10)}$), Lemma \ref{lem:suma} (with $f(a)=\tau(a)$), \eqref{apz} and \eqref{plQ} that the number of integers $n\in (x/2,x]$ counted
for a given $z$ and $l$ is at most $T_1(z,l)+T_2(z,l)+T_3(z,l)$, where
\begin{align*}
T_1(z,l) &= \frac{xe^{-l}}{\log^2 z} \sum_{a\in \sP(z)} \frac{L_{\ge \beta r}(a)}{a}, \qquad \beta= \frac{1}{160l^2}, \\
T_2(z,l) &=  \frac{x}{\log^2 z} \sum_{a\in \sP(z^{\exp(-l/10)})} \frac{L_{\ge \beta r}(a)}{a}, \\
T_3(z,l) &= \frac{x}{(\log z)^{100}} \sum_{a\in\sP(z)}\frac{\tau(a)}{a} \ll \frac{x}{(\log z)^{98}}.
\end{align*}
Since $r\le \log z$ and $x\ge Y/2$,
\[
\sum_{l=1}^r T_3(z,l) \ll \frac{x}{(\log z)^{97}} \ll \frac{x}{(\log Y)^{97}}.
\]
Lemma \ref{lem:suml1} implies that
\begin{align*}
\sum_{l=1}^r T_1(z,l) &\ll \frac{x\log(2r)}{(\log z)^{1-E(\vartheta(r,z))}(\log\log z)^{3/2}},\\
\sum_{l=1}^r T_2(z,l) &\ll \frac{x\log(2r)}{(\log z)^{1-E(\vartheta(r,z))}(\log\log z)^{3/2}}.
\end{align*}
The function $(E(x)-1)/x$ is increasing for $0<x\le \log 4-1$ and decreasing thereafter,
thus the right sides above are decreasing in $z$ upon taking $x=\vartheta(r,z)$; 
incidentally, this is one place in our proof where 
we cannot take $r$ larger than $(\log Y)^{\log 4-1}$.
Recalling that $z\ge Y/2$, this completes the proof.
\end{proof}

\begin{cor}\label{cor:H-tilde}
Fix $\delta'''>0$.
For $Y$ large enough in terms of $\delta'''$, where $Y=\min(y,x/(2y))$, and 
$r  \leq (\log Y)^{\log 4-1-\delta'''},$ we have 
\begin{equation}
\tH_{\ge r}(x,y,2y) \ll_{\delta'''} \frac{x \log(2r)}{(\log Y)^{1-E(\vartheta(r,Y))}(\log\log Y)^{3/2}}.
\end{equation}
\end{cor}

\begin{proof}
For $j\ge 0$ with $2^j\le \log Y$ let 
\[
x_j = \frac{x}{2^j}, \qquad Y_j = \min\bigg(y,\frac{x}{2^{j+1}y}\bigg), \qquad
\theta_j = \vartheta(r,Y_j).
\]
Then
\begin{equation}\label{H-xj}
\tH_{\ge r}(x,y,2y)\le \ssum{j\ge 0 \\ 2^j\le \log Y} \Big( \tH_{\ge r}(x_j,y,2y)-\tH_{\ge r}(x_{j+1},y,2y)\Big) + O\pfrac{x}{\log Y}.
\end{equation}
The error term is acceptable since $E(\vartheta(r,Y))\ge E(\log 4-1)=2-\log 4>0.5$.
Let $\delta''=\frac12 \delta'''$.
We apply Lemma \ref{lem:H-tilde-short-int} with $x=x_j$, noting that $Y_j\gg Y/\log Y$
and hence $\theta_j \le \log 4 -1+\delta''$ if $Y$ is large enough in terms of $\delta'''$.
Also, since $Y_j \le Y$, $\theta_j \ge \vartheta(r,Y)$ and hence $E(\theta_j)\le E(\vartheta(r,Y))$.
It follows from Lemma \ref{lem:H-tilde-short-int} that
\[
 \tH_{\ge r}(x_j,y,2y)-\tH_{\ge r}(x_{j+1},y,2y) \ll_{\delta'''} \frac{2^{-j} x \log(2r)}{(\log Y)^{1-E(\vartheta(r,Y))}(\log\log Y)^{3/2}}.
\]
Summing on $j$ completes the proof.
\end{proof}

Now we are ready to prove the upper bound in Theorem \ref{th1.1}. Without loss of generality, we may think that $Y$ is large enough in terms of $\delta$, since the upper bound in
Theorem \ref{th1.1} is trivial for bounded $Y$.
 For $n$  counted in $H_{\ge r}(x,y,2y)$, we write $n=n'n''$, where $n'$ is squarefree, $n''$ is square-full and $\gcd(n',n'')=1$.  The number of $n$ with $n''>\log^2 Y$ is $O(x/\log Y)$.
 Let $1<n''\le \log^2 Y$, so that $n''\geq 4$. Every divisor $d_i$ of $n$ can be written as $d_i=d_i'd_i''$, where $d_i'|n'$ and $d_i''|n''$. Each $d_i'$ belongs to the interval 
 $(y/n'',2y]$ and thus to some interval $(y2^{-j+1},y2^{-j+2}]$ with $1\leq j\leq (\log n'')/(\log2)+2 \le 3\log n''$.  By the pigeonhole principle, for some such $j$ there 
 are at least $r/(3\log n'')$ indices $i$ for which $d_i'\in (y2^{-j+1},y2^{-j+2}]$.
For any $d|n'$, at most $\tau(n'')$ of these $d_i'$ equal $d$.
We obtain
\[
H_{\ge r}(x,y,2y) \leq \ssum{4\le n'' \le \log^2 Y \\ n''\text{ square-full}}\; \sum_{2^j\leq 4n''}\tH_{\tr}\left(\frac{x}{n''},\frac{y}{2^{j-1}},\frac{y}{2^{j-2}}\right) +\tH_r(x,y,2y)+O\pfrac{x}{\log Y}.,
\]
where
\[
\tr = \tr(r,n'') = \cl{\frac{r}{3\tau(n'') \log n''}}.
\]
Fix $n''$ and $j$ and let
\[
\tY= \min\bigg( \frac{y}{2^{j-1}}, \frac{x/n''}{2y/2^{j-1}}\bigg), \qquad \ttheta=\vartheta(\tr,\tY).
\]
Since $2^j\le 4n'' \le 4\log^2 Y$,
\[
\frac{Y}{2\log^2 Y} \le \tY \le 2Y.
\]
  If $Y$ is large enough then
 $\ttheta \le \theta +\frac12 \delta \le \log 4 -1 - \frac12 \delta$.  Applying
 Corollary \ref{cor:H-tilde} with $\delta'''=\frac12 \delta$, we have
 \begin{equation}\label{Hrn''j}
 \tH_{\tr}\bigg(\frac{x}{n''},\frac{y}{2^{j-1}},\frac{y}{2^{j-2}}\bigg)
 \ll_{\delta} \frac{x\log(2r)}{n'' (\log Y)^{1-E(\ttheta)}(\log\log Y)^{3/2}}. 
 \end{equation}
 Now 
\[
 \ttheta = \frac{\log r + O(\log(3\tau(n'')\log n''))}{\log\log Y+O(1)} = 
 \theta + O\pfrac{1+\log (\tau(n'')\log n'')}{\log\log Y},
\]
hence
\[
(\log Y)^{1+E(\ttheta)} = (\log Y)^{1+E(\theta)} (2\tau(n'')\log n'')^{O(1)}.
\]
 Inserting this into \eqref{Hrn''j}, using that
 \[
 \sum_{n''} \frac{(2\tau(n'')\log n'')^{O(1)}}{n''} \ll 1,
 \]
 and applying Corollary \ref{cor:H-tilde} for the term
 $\tH_{\ge r}(x,y,2y)$,  the proof of the upper bound in Theorem \ref{th1.1} is complete.

\medskip

It remains to prove Proposition \ref{prop1.2}. This will be achieved in the next section.

\medskip

\section{Proof of Proposition \ref{prop1.2}: upper bound} \label{Sec3}

Recall that our goal is to show that 
\begin{equation}\label{3.1}
\sum_{a\in \sP(z)}\frac{L_{\ge t}(a)}{a} \ll \frac{(\log (2t))(\log z)^{1+E(\vartheta(t,z))}}
{(\log\log z)^{3/2}}
\end{equation}
uniformly for $z\geq 3$ and $1\leq t \leq (\log z)^{\log 4-1-\delta'}$.  
If $L_{\ge t}(a)>0$ then $\tau(a)\ge t$ and thus $\omega(a)\ge \frac{\log t}{\log 2}$.
Thus
\[
\sum_{a\in \sP(z)}\frac{L_{\ge t}(a)}{a}=\sum_{k\geq \frac{\log t}{\log 2}}T_k,
\]
where 
\begin{equation*}
T_k=\sum_{\substack{a\in \sP(z)\\ \omega(a)=k}}\frac{L_{\ge t}(a)}{a}.	
\end{equation*}

For $a$ with $\omega(a)=k$, Lemma \ref{lem:Lr}(a) gives the crude bound 
\begin{equation}\label{Lr-crude}
L_{\ge t}(a)\le \min\{2^{k+2}/t, \log (2a)\}. 
\end{equation}
Let
\begin{equation}\label{k0v}
k_0=\left\lfloor\frac{\log\log z+\log t}{\log2}\right\rfloor, \qquad
v=\frac{\log\log z}{\log 2}.
\end{equation}
Since most $a$ with $P^+(a)\leq z$ obey $\log a\ll \log z$, one can expect that the critical value of $k$ is $k=k_0$.
The next lemma allows us to concentrate on the values of $k$ which are close to $k_0$.

\begin{lem}\label{lem:Tk tails}
Fix $\delta'>0$.  Uniformly for $1\le t\le (\log z)^{\log 4-1-\delta'}$ and $k_1\le k_0$,
\[
\sum_{k\le k_1} T_k \ll_{\delta'} \frac{(2\log\log z)^{k_1}}{k_1! t}.
\]
Uniformly for $t\ge 1$ and $k_2\ge k_0$, we have
\[
\sum_{k\ge k_2} T_k \ll \frac{(\log z)(\log\log z)^{k_2-1}}{(k_2-1)!}.
\]
\end{lem}

\begin{proof}
For each $k\le k_0$, \eqref{Lr-crude} gives
\begin{equation}\label{Tk-smallk}
T_k\leq \frac{2^{k+2}}{k!t}\left(\sum_{p\leq z}\frac1p\right)^k \ll \frac{(2\log\log z+O(1))^k}{k!t}.
\end{equation}
Sum over $k\le k_1$ and note that $k_1\le k_0\ll \log\log z$ implies
$(2\log\log z+O(1))^k\ll(2\log\log z)^k$ for $k\le k_0$.
Finally, note that $k_1\le k_0< (2-\delta'/\log 2)\log\log z+O(1)$ and the first bound follows.

For $k>k_0$, \eqref{Lr-crude} gives
\begin{align*}
T_k &\ll \ssum{a\in \sP(z)\\ \omega(a)=k} \frac{\log a}{a} = 
\ssum{a\in \sP(z)\\ \omega(a)=k} \frac{1}{a} \sum_{p|a} \log p \\
&\le \sum_{p\le z} \frac{\log p}{p} \ssum{b\in \sP(z)\\ \omega(b)=k-1} \frac{1}{b}
\le \bigg(\sum_{p\le z} \frac{\log p}{p}\bigg) \, \frac{1}{(k-1)!} \bigg( \sum_{p\le z} \frac{1}{p}\bigg)^{k-1}\\
&\ll \frac{(\log z)(\log\log z+O(1))^{k-1}}{(k-1)!}.
\end{align*}
Since $k_0 > 1.3\log\log z$, the sum over $k \ge k_2$ of the right side
 is big-$O$ of the single term with $k=k_2$.
\end{proof}

\textbf{Proof of Proposition \ref{prop1.2}(a).}
Assume first that $(\log z)^{1/100} < t < (\log z)^{\log 4-1-\delta'}$.
Take $k_1=k_0$ and $k_2=k_0+1$ in Lemma \ref{lem:Tk tails}, noting that the right side in both estimates
is $\order (\log z)^{1+E(\vartheta(t,z))}(\log\log z)^{-1/2}$, matching the bound in 
Proposition \ref{prop1.2}(a) since $\log(2t) \gg_{\delta'} \log\log z$.

\medskip

The case $t<(\log z)^{1/100}$ is more delicate.
Here $1.4 \log\log z \le k_0 \le 1.5\log\log z$.  Thus, taking 
$k_1 = 0.9 k_0$ and $k_2 = 1.1 k_0$, we find that
\begin{equation}\label{eq:Tk-tails}
\sum_k T_k \le \sum_{0.9k_0 \le k\le 1.1k_0} T_k + O\pfrac{(\log (2t))(\log z)^{1+E(\vartheta(t,z))}}{(\log\log z)^2}.
\end{equation}

As in \cite[Section 4]{F1}, partition the primes as follows.
Let $\lambda=1.9$ and, for $j\geq1$, inductively define $\lambda$ as the largest
prime so that
\begin{equation}\label{lamj-def}
\sum_{\lambda_{j-1}<p\leq \lambda_j}\frac1p \leq \log 2.
\end{equation}
Lemma 4.6 of \cite{F1} implies
\begin{equation}\label{lamj-approx}
\log \lambda_i \asymp 2^i.	
\end{equation}
Let $D_j=\{p: \lambda_{j-1}<p\leq \lambda_j\}$. Now we estimate each sum $T_k$ in terms of a multivariate integral. 

\begin{lem}\label{lem:Tk-medium-k}
Fix $B\ge 2$.  Uniformly for $k\le Bv$, $z\ge 3$, $t\ge 1$,
\[
T_k \ll_B (2\log\log z)^kU_k(v,t),
\]
where 
\begin{equation*}\label{Ukv}
U_k(v,t)=\idotsint\limits_{0\leq \xi_1\leq \cdots \leq \xi_k\leq1}
\min\left(\frac1t,\min_{1\leq j\leq k}2^{-j}(2^{v\xi_1}+\cdots+2^{v\xi_j})\right) d\bxi
\end{equation*}
\end{lem}

\begin{proof} For a number $a=p_1\cdots p_k$ with $p_1<\cdots<p_k\leq z$, we define $m_i$ by $p_i\in D_{m_i}$.  By \eqref{lamj-approx}, $\log p_i \ll 2^{m_i}$ and
\[
m_k\leq \max\{j: \lambda_{j-1}\leq z\} =: v',
\]
where $v'=v+O(1)$. By Lemma \ref{lem:Lr}(c),
\begin{align*}
L_{\ge t}(a) &\ll 2^{k}\min\left(\frac1t,\min_{1\leq j\leq k}2^{-j}(2^{m_1}+\cdots+2^{m_j}+\log 2)\right)\\
&\ll 2^k\min\left(\frac1t,\min_{1\leq j\leq k}2^{-j}(2^{m_1}+\cdots+2^{m_j})\right)  = 2^k F(\bm),
\end{align*}
where, for a real vector $\mathbf{s}=(s_1,\ldots,s_k)$ with $0\leq s_1\leq\cdots\leq s_k$, 
\[
F(\mathbf{s})=\min\left(\frac1t,\min_{1\leq j\leq k}2^{-j}(2^{s_1}+\cdots+2^{s_j})\right). 
\]
For a fixed $\bm$, let $b_i=b_i(\bm) = \#\{u : p_u \in D_{m_i}\}$ and also $\beta_0=0$ and $\beta_i=b_1 + \cdots + b_i$ for $i\geq1$. In particular, $\beta_{v'}=b_1+\cdots+b_{v'}=k$.  Then 
\begin{align*}
T_k &\ll 2^k\sum_{ 1\leq m_1\leq \cdots\leq m_k\leq v'}F(\bm)
\prod_{j=1}^{v'}\bigg(\sum_{\substack{p_{\beta_{j-1}+1}<\cdots<p_{\beta_j}\\
p_{\beta_{j-1}+1},\ldots,p_{\beta_j} \in D_j}}
\frac{1}{p_{\beta_{j-1}}\cdots p_{\beta_j}}\bigg) 
\\ &\leq 2^k\sum_{1\leq m_1\leq \cdots \leq m_k\leq v'}F(\bm)
\prod_{j=1}^{v'} \frac{(\log 2)^{b_j}}{b_j!}
= (2 \log 2)^k\sum_{1\leq m_1\leq \cdots \leq m_k\leq v'}
\frac{F(\bm)}{b_1!\cdots b_{v'}!}.
\end{align*}
If $m_i-1 \leq s_i \leq m_i$ for each $i$, then $F(\bm) \leq 2F(\mathbf{s})$. Therefore,
\[
\sum_{1\leq m_1\leq \cdots \leq m_k\leq v'}
\frac{F(\bm)}{b_1!\cdots b_{v'}!}\leq 2\idotsint\limits_{0\leq s_1\leq \cdots \leq s_k\leq v'}F(\mathbf{s})d\mathbf{s}
\]
Making the change of variables $s_i = v'\xi_i$ for each $i$, 
\[
\idotsint\limits_{0\leq s_1\leq \cdots \leq s_k\leq v'}F(\mathbf{s})d\mathbf{s} =
(v')^k\idotsint\limits_{0\leq \xi_1\leq \cdots \leq \xi_k\leq1}
\min\left(\frac1t,\min_{1\leq j\leq k}2^{-j}(2^{v'\xi_1}+\cdots+2^{v'\xi_j})\right) d\bxi.
\]
Recalling that $v' = v+ O(1)$, we conclude that
\[
T_k \ll (2v'\log 2)^k U_k(v,t),
\]
and the claim follows from 
\[
(2v'\log 2)^k = (2v\log 2 + O(1))^k \ll_B (2\log\log z)^k
\]
since $k\leq Bv$.
\end{proof}

The integrals $U_k(v,t)$ can be interpreted as follows: let $Y_1,\ldots,Y_k$ be independent samples of a random variable which has uniform distribution on $[0,1]$, and let $\xi_1\leq \xi_2\leq \cdots \leq \xi_k$ be their non-decreasing permutation. For each $k$, define the random variable
$$
X_k=\min\left(\frac1t,\min_{1\leq j\leq k}2^{-j}(2^{v\xi_1}+\cdots+2^{v\xi_j})\right). 
$$
Then (with $\mathbb{E}$ denoting expectation)
\begin{equation*}
k!U_k(v,t) = \mathbb{E} X_k.
\end{equation*}
Roughly, we expect that, most of the time,
\begin{equation*}
2^{v\xi_1}+\cdots+2^{v\xi_j} \ll 2^{v\xi_j}.
\end{equation*}
Assuming this, one can deduce (\ref{3.1}) from the following result, which is an immediate consequence of \cite[Theorem 1]{F2}. 

\begin{lem}\label{lem3.3}
Uniformly in $u\geq 0$ and $u+v \geq k$, we have ($\mathbb{P}$ denotes probability)
\[
\mathbb{P}\left(v\xi_i-i \geq -u \mbox{ for all } i\right) \ll \frac{(u+1)(u+v-k+1)}{k}.
\]	
\end{lem}

We omit the details of the mentioned conditional claim and just proceed further with the proof. To estimate $U_k(v)$, we bound the volume of the set
\[
\T(k,v,m)=\{\bxi \in \RR^k: 0\leq \xi_1\leq \cdots \leq \xi_k\leq 1: 2^{v\xi_1}+\cdots+2^{v\xi_j}\geq 2^{j-m} \, (1\leq j\leq k)\}. 
\]
Inserting Lemma \ref{lem3.3} into the proof of Lemma 13.2 of \cite{F1} (with no other changes to the proof), we get
\begin{lem}\label{lem3.4}
Suppose that $k,v,m$ are integers with $1\leq k\leq 10v$ and $m\geq0$. Set $b=k-v$. Then
\[
\Vol(\T(k,v,m)) \ll \frac{Y_m}{2^{2^{b-m}}(k+1)!}, \quad 
Y_m=\begin{cases}
b, & \mbox{ if } \, b\geq m+5, \\
(m+5-b)(m+1), & \mbox{ if } \, b\leq m+4.	
\end{cases} 
\]
\end{lem}

Now we are ready to finish the proof of (\ref{3.1}). Let 
\[
m_0=\left\lfloor \frac{\log t}{\log 2}\right\rfloor .
\]
In the integral defining $U_k(v)$, the minimum is in $[2^{-m},2^{1-m}]$ for some $m\ge m_0$.
Hence, for $0.5v\leq k\leq 3v$, we have by Lemma \ref{lem3.4},
\[
U_k(v,t) \ll \sum_{m\geq m_0}2^{-m}\Vol(\T(k,v,m)) \ll \frac{1}{(k+1)!}\sum_{m\geq m_0}\frac{2^{-m}Y_m}{2^{2^{b-m}}}.
\]
Further,
\[
\sum_{m\geq m_0}\frac{2^{-m}Y_m}{2^{2^{b-m}}} \ll \sum_{m_0\leq m\leq b-5}\frac{b}{2^{m+2^{b-m}}}+\sum_{m\geq \max\{m_0,b-4\}}\frac{(m+5-b)(m+1)}{2^m}
\]
If $k\ge k_0+5$ then  $b\geq m_0+4$ and both sums above are $\ll b2^{-b}$. 
If $k\le k_0+4$, then $b-5<m_0$ and the first sum is empty, while the second sum is 
 \[
 \ll 2^{-m_0}(m_0+1)(m_0+5-b) \ll \frac{\log(2t)}{t} (k_0-k+5).
\]
  Putting everything together and noting that $2^{-m_0} \order 1/t$, we get
\begin{equation*}
(k+1)! U_k(v,t) \ll \begin{cases}
\frac{(k-k_0+m_0)}{t\cdot 2^{k-k_0}}, & \mbox{ if } \, k\geq k_0+5, \\
\\
\frac{(\log (2t))(k_0-k+5)}{t}, & \mbox{ if } \, k\le k_0+4.
\end{cases}
\end{equation*}
Hence, by Lemma \ref{lem:Tk-medium-k} and writing $(k+1)!\order (\log\log z) k!$,
\begin{align*}
\sum_{0.9k_0\leq k\leq 1.1k_0}T_k 
\ll \frac{\log (2t)}{t \log\log z}&\;\sum_{0.9k_0\leq k\le k_0+4}\frac{(2\log\log z)^k (k_0-k+5)}{k!}+\\
&\qquad+\frac{2^{k_0}}{t \log\log z}\;\sum_{k_0+5\leq k  \leq 1.1k_0}
\frac{(\log\log z)^k (k-k_0+m_0)}{k!}.	
\end{align*}
Since $k_0\le (2-\delta'/\log 2)\log\log z$ by hypothesis,
both sums are dominated by terms with $k=k_0+O(1)$, hence
\begin{equation*}\label{3.6}
\sum_{0.9k_0\leq k\leq 1.1k_0}T_k \ll \frac{\log(2t)}{t \log\log z}\frac{(2\log\log z)^{k_0}}{k_0!} \asymp \frac{(\log(2t)) (\log z)^{1+E(\vartheta(t,z))}}{(\log\log z)^{3/2}}.
\end{equation*}
Together with \eqref{eq:Tk-tails},
this concludes the proof of Proposition \ref{prop1.2}(a).

\section{Proof of Proposition \ref{prop1.2}: lower bound}

In this section we establish the lower bound in Proposition \ref{prop1.2} (b), namely
\begin{equation}\label{4.1}
\sum_{a\leq z}\frac{L_{\ge t}(a)}{a} \gg \frac{(\log z)^{1+E(\vartheta(t,z))}\log (2t)}{(\log\log z)^{3/2}}
\end{equation}
uniformly in large enough $z$ and $1\leq t\leq (\log z)^{10}$.
Let 
\[
k_0=\fl{\frac{\log(t\log z)}{\log 2}}+4, \qquad \sigma=\frac{1}{2Dt},
\]
 where $D$ is a large constant to be chosen later.
We will show that the contribution from $a$ with $\omega(a)=k_0$ gives the above lower bound. For $a\in\N$, we define
\begin{equation*}
W(a;\sigma)=\#\{(d_1,d_2): d_1|a, d_2|a, |\log(d_1/d_2)|\leq \sigma\}.	
\end{equation*}
In particular, the diagonal terms $d_1=d_2$ give $W(a;\sigma)\ge \tau(a)$.

\begin{lem}\label{lem:Lt-lower}
Let $a\leq z$ be such that $\omega(a)=k_0$, $a$ is squarefree, and $W(a,\sigma)\leq \frac32\tau(a)$. Then	
\[
L_{\ge t}(a) \geq \frac{0.3}{D}\log z.
\]
\end{lem}

\begin{proof}
Let $S$ be the set of \emph{isolated} divisors $d|a$, in the sense that $|\log(d/d')|>\sigma$
for all divisors $d'|a$ with $d\ne d'$.  Evidently,
\[
W(a;\sigma) \ge \tau(a) + (\tau(a)-|S|),
\]
whence $|S| \ge \frac12 \tau(a)$.

Now for $k\geq0$ we set $I_k=[e^{k/2},e^{(k+1)/2})$ and write
\begin{equation}\label{4.2}
|S|=\sum_{k: |S\cap I_k|<t} |S\cap I_k| + \sum_{k: |S\cap I_k|\geq t}\, |(S\cap I_k)|,
\end{equation}
where $-2 \le k\le 2\log z$.  In particular, the number of $k$ under consideration is 
at most $2.1\log z$.
Note that $|S\cap I_k|\leq 0.5\sigma^{-1}=Dt$ for any $k$, because if $d',d''$ are elements of $S$ with $d'>d''$, then $d'>d''e^{1/\sigma}$. Since $\tau(a) =  2^{k_0}\geq 8t\log z$, we get from (\ref{4.2}) 
\[
4t\log z \le \tfrac12 \tau(a) \leq t(2.1\log z)+\big| \{k: |S\cap I_k|\geq t\}\big|\cdot Dt,
\]
so that	at least $\frac{1.9}{D}\log z$ of the intervals $I_k$ each 
contain at least $t$ elements of $S$. For each $k$, if $x\in [k/2+\frac12-\log2, k/2)$, then the interval $(e^x,2e^x]$ contains $I_k$. Thus,
\[
L_{\ge t}(a)\geq (\log2-\tfrac12) \frac{1.9 \log z}{D} \ge \frac{0.3}{D}\log z.\qedhere
\]
\end{proof}

Lemma \ref{lem:Lt-lower} gives
\begin{equation}\label{4.3}
\sum_{a\leq z}\frac{L_{\ge t}(a)}{a} \geq \frac{0.3}{D}(\log z)\sum_{\substack{a\leq z, \,\, \omega(a)=k_0, \\ W(a,\sigma)\leq \frac32 \tau(a)}} \frac1a \geq \frac{0.3\log z}{2^{k_0}D}\sum_{\substack{a\leq z\\\omega(a)=k_0}}\frac{3\tau(a)-2W(a;\sigma)}{a}.,	
\end{equation}
where we used $3\tau(a)-2W(a;\sigma) \le \tau(a) = 2^{k_0}$ in the last inequality.

We now quote some results from \cite[\S 4]{F1}.  Here $c_5$ is a constant 
appearing in Lemmas 4.7 and 4.8 of \cite{F1},
and $M$ is a sufficiently large even integer which is needed for
 \cite[Lemmas 4.8, 4.9]{F1}.
 If $z$ is large enough then $0<\sigma<2^{-2M-1}\log z$.
Hence, by \cite[Lemmas 4.8]{F1} (with $\alpha=1$), the sum on the right side of \eqref{4.3} is 
\[
\gg (2v\log2)^{k_0}\Vol(Y_{k_0}(s,v)),
\]
where 
\begin{align*} 
s=\frac{\log(2Dt)}{\log 2}+M-c_5-10, \quad
v=\left\lfloor \frac{\log\log z}{\log2}\right\rfloor-2M, 
\end{align*}
and $Y_{k_0}(s,v)$ is the set of $\bxi = (\xi_1,\ldots, \xi_{k_0})\in \RR^{k_0}$ satisfying
\begin{align*}
\tag{i} & 0\leq \xi_1 \leq \cdots \leq \xi_{k_0} \leq 1; \\ 
\tag{ii} & \xi_{M+i^2}>i/v \mbox{ and } \xi_{k_0+1-(M+i^2)}<1-i/v  \mbox{ for } 1\leq i \leq\sqrt{k_0-M}; \\
\tag{iii} &\sum_{j=1}^{k_0}2^{j-v\xi_j} \leq 2^s.
\end{align*}

\begin{lem}\label{lem4.9-var}
\cite[Lemma 4.9]{F1}.
For $v\ge 1$, $10M\le k\le 100(v-1)$, $s\ge M/2+1$ and $0\le k-v\le s-M/3-1$, we have
\[
\Vol \big( Y_k(s,v) \big) \gg \frac{k-v+1}{(k+1)!}.
\]
\end{lem}
%
% this is for a version with 100 replaced by arbitrary B.  Don't need
%\begin{proof}
%Lemma 4.9 of \cite{F1} provides the case $B=98$ and there is no essential difference in
%replacing $100$ with an arbitrary $B+2$. The argument is on pages 423--425 of \cite[Section 12]{F1}, and inequality and Lemma numbers refer to displays in that paper.
%In (12.5) and (12.6), simply replace 100 with $B+2$; here the implied constants 
%in the $\ll$ symbols are absolute.
%Three lines before (12.7), argue as follows: since $k\le (B+2)(v-1)$,
%$v-1-k\ge (-1+1/(B+2))k$ and apply Lemma 12.1 with $\eps=1/(B+2)$ instead of $\eps=1/100$.
%As Lemma 12.1 has a factor $e^{5/\eps}$, this induces a multiplicative factor in (12.7) 
%which depends on $B$.  Replacing 100 by $B$ in (12.9) forces the implied constant
%in (12.9) to depend on $B$.  At the top of page 425, the constants
%$K_1$ and $K_2$ depend on $B$, forcing $M$ to be sufficiently large in terms of $B$
%for the required estimate to hold.
%\end{proof}

Let $k=k_0$. The first conditions $v\ge 1$ and $10M\le k\le 100(v-1)$ are all met if $z$ is large enough in terms of $M$, and
the conditions $s\ge M/2+1$ and $0\le k-v\le s-M/3-1$  are met if
 $D$ is large enough (note that we may choose $D$ dependent on $M$).
Combining everything together, we find
\[
\sum_{a\leq z}\frac{L_{\ge t}(a)}{a}  \gg \frac{(\log z)\log (2t)}{\log\log z}\frac{(v\log 2)^{k_0}}{k_0!}
\]
Since $1\le t\le (\log z)^{10}$, $k_0\ll \log\log z$ and hence
$(v\log 2)^{k_0}\order (\log\log z)^{k_0}$.  Now $k_0=\frac{1+\theta}{\log 2}\log\log z+O(1)$,
so by Stirling's formula, the right side is
\[
\gg \frac{(\log z)\log(2t)}{(\log\log z)^{3/2}} \pfrac{e\log 2}{1+\theta}^{\frac{1+\theta}{\log 2}\log\log z} = \frac{(\log z)^{1+E(\vartheta(t,z))}\log(2t)}{(\log\log z)^{3/2}},
\]
proving claim (\ref{4.1}).

\bigskip

\textbf{Acknowledgements.}   The first autho was supported in part by
National Science Foundation grant DMS-2301264, a Simons Foundation Fellowship
and a Simons Foundation Travel Support grant.
The third author is partially supported by the National Science Foundation (DMS-2418328) and the Simons Foundation (MPS-TSM-00007959).

\end{document}